# *Double Fuzzy Complex EE Transform to Solve Partial Volterra Fuzzy Integro-Differential Equations*

*Jinan A. Jasim*

*Department of Mathematics*

*College of Education*

*Mustansiriyah University*

*Baghdad, Iraq*

jinanadel78@uomustansiriyah.edu.iq

*Alan Jalal Abdulqader*

*Department of Mathematics*

*College of Education*

*Mustansiriyah University*

*Baghdad, Iraq*

alanjalal515@uomustansiriyah.edu.iq

*Emad A. Kuffi*

*Department of Materials*

*College of Engineering*

*Al-Qadisiyah University*

*Al-Qadisiyah, Iraq*

emad.abbas@qu.edu.iq

***Abstract***: In this paper, the double fuzzy complex EE transform was applied to get the solution to partial Volterra fuzzy integro-differential equations with convolution kernel under H-differentiability. This work presents important results to this transform for double fuzzy convolution and fuzzy partial derivatives of the n-th order that help us to solve previously mentioned equations . a given example shows the method of applying DFCEET for solving Volterra Fuzzy Integro-Differential Equations.

***Introduction:*** Modeling of different physical systems gives us different differential, integral and integro-differential equations. We are not always sure that the models obtained are perfect. The fuzzy set theory is one of the most popular theories for describing this situation. The fuzzy logic is introduced with the proposal of fuzzy set theory by Zadeh [1] and is applied when the observational parameters are imprecise or unclear. The neutrosophic logic is considered as the extension of the fuzzy logic and the measure of indeterminacy is added to the measures of truthiness and falseness. The theory of neutrosophic statistics can be applied when to observation indeterminate, imprecise, vague, and incomplete parameters[2-5]. In last years, many mathematicians have studied the solution of fuzzy differential equations, fuzzy integral equations, and fuzzy integro-differential equations, which play a key role in engineering. These equations in a fuzzy setting are a natural way to model the ambiguity of dynamic systems in different scientific fields [6-15].

The integral transforms be very useful in solving partial differential equations. They convert the original function to a function that is simpler to solve. The Fourier transform is the precursor of the integral transforms. This transform is used to express functions in a finite interval. Similar integral transforms are Laplace, Sumudu, Al-Tememe, Complex Al-Tememe, Mahgoub, Mohanad, SEE, complex SEE and complex EE transforms [16-24]. Many researches discussed solving of fuzzy integral equations in all kinds by using fuzzy integral transforms or double fuzzy integral transforms for equations with two variables [25-27].

In the present paper we developed the double complex EE transform [28] to fuzzy double complex EE (DFCEET) transform and solve the fuzzy partial Volterra integro-differential convolution integral equation of two variable under Hukuhara differentiability using (DFCEET) .The main difficulties overcome in solving the this problem are related to the application of the DFCEET for fuzzy partial H-derivative of the n-th order. So, we obtain a new results on the complex EE transform for fuzzy partial H-derivative of the n-th order . After, the studied equation we convert to a nonlinear system of partial Volterra integro-differential equations in a crisp case. To be find the lower and upper functions of the solution we use DFCEET and we convert this system to system of algebraic equations. The paper is introduced firstly some definitions and results of fuzzy numbers, fuzzy functions and fuzzy partial derivative of the n-th order is given. After that the definition of DFCEET is recalled, double fuzzy convolution theorem is stated. New results on DFCEET for fuzzy partial derivative of the n-th order are proposed. The DFCEET is applied to fuzzy partial convolution Volterra fuzzy integro-differential equation to construct the general algorithm who steps are illustrated by a given example. Finally, conclusions are drawn.

***Definition (1)[27]:*** A fuzzy number is a function $\tilde{v}: \mathbb{R} \to [0, 1]$ satisfying the following properties:

i. $\tilde{v}$ is an upper semi-continuous on $\mathbb{R}$.
ii. $\tilde{v}(x) = 0$ outside of some interval $[c, d]$.
iii. There are the real numbers $a$ and $b$ with $c \leq a \leq b \leq d$, such that $\tilde{v}$ is increasing on $[c, a]$, decreasing on $[b, d]$ and $\tilde{v}(x) = 1$ for each $x \in [a, b]$.
iv. $\tilde{v}$ is a fuzzy convex set that is $\tilde{v}(\alpha x + (1 - \alpha)y) \geq \min\{\tilde{v}(x), \tilde{v}(y)\}$ for all $x, y \in \mathbb{R}$ and $\alpha \in [0, 1]$.
v. The set of all fuzzy numbers is denoted by $F^1$ and $D = \mathbb{R}_+ \times \mathbb{R}_+$. Any real number $r \in \mathbb{R}$ can be interpreted as a fuzzy number $\tilde{r} = X_{[r]}$, therefore $\mathbb{R} \subset F^1$. $\mathbb{R}_f$ denote the set of all fuzzy numbers on $\mathbb{R}$.
vi. $\underline{v}(\alpha)$ is a bounded non-decreasing left continuous function in $(0, 1]$, and right continuous at $0$.
vii. $\bar{v}(\alpha)$ is a bounded non-increasing left continuous function in $(0, 1]$, and right continuous at $0$.
viii. $\underline{v}(\alpha) \leq \bar{v}(\alpha), \ \alpha \in [0, 1]$.

***Definition (2) [29]***: for arbitrary $\tilde{v}, \tilde{z} \in F^1$ such that $\tilde{v} = \left(\underline{v}(\alpha), \bar{v}(\alpha)\right), \tilde{z} = \left(\underline{z}(\alpha), \bar{z}(\alpha)\right), \ \alpha \in [0, 1]$ we define:

1. Addition $\quad \tilde{v} \oplus \tilde{z} = \left(\underline{v}(\alpha) + \underline{z}(\alpha), \bar{v}(\alpha) + \bar{z}(\alpha)\right)$.

2. Subtraction $\tilde{v} \ominus \tilde{z} = \left(\underline{v}(\alpha) - \underline{z}(\alpha), \overline{v}(\alpha) - \overline{z}(\alpha)\right)$.
3. Multiplication $\tilde{v} \odot \tilde{z} =$
$\left(\min\{\underline{v}(\alpha)\overline{z}(\alpha), \underline{v}(\alpha)\underline{z}(\alpha), \overline{v}(\alpha)\overline{z}(\alpha), \overline{v}(\alpha)\underline{z}(\alpha)\}, \max\{\underline{v}(\alpha)\overline{z}(\alpha), \underline{v}(\alpha)\underline{z}(\alpha), \overline{v}(\alpha)\overline{z}(\alpha), \overline{v}(\alpha)\underline{z}(\alpha)\}\right)$.

***Definition (3) [2 7]:*** let $\tilde{w}: D \to F^1$ be a fuzzy-number-valued function. Then $w$ is said to be continuous at $(s_0, t_0) \in D$ if for each $\epsilon > 0$ there is $\delta > 0$ such that $d(f(s,t), f(s_0, t_0)) < \epsilon$ whenever $|s - s_0| + |t - t_0| < \delta$. If $w$ be continuous for each $(s, t) \in D$ then we say that $w$ be continuous on $D$.

Let $R > 0$. Denote $D_R = D \cap \overline{U}(0, R)$, where

$$\overline{U}(0, R) = \{(x, y): x^2 + y^2 \leq R^2\}$$

is the closed circle with radius $R$.

***Theorem (1) [2 7]:*** let $\tilde{w}: D \to F^1$ be fuzzy-valued function with parametric form $\left(\underline{w}(x, y, \alpha), \overline{w}(x, y, \alpha)\right)$ for all $\alpha \in [0, 1]$:

1. The functions $\underline{w}(x, y, \alpha)$ and $\overline{w}(x, y, \alpha)$ are Riemann-integrable on $D_\mathbb{R}$.
2. There are constants $\underline{M}(\alpha) > 0$ and $\overline{M}(\alpha) > 0$, such that:

$$\iint_{D_R} [\underline{w}(x, y, \alpha)] \, dx \, dy \leq \underline{M}(\alpha)$$

$$\iint_{D_R} [\overline{w}(x, y, \alpha)] \, dx \, dy \leq \overline{M}(\alpha), \quad \forall R > 0,$$

Then the function $\tilde{w}(x, y)$ is an improper fuzzy Riemann-integrable on $D$ and

$$(FR)\int_0^\infty (FR)\int_0^\infty \tilde{w}(x,y)\,dx\,dy = \left(\int_0^\infty\int_0^\infty \underline{w}(x,y,\alpha)\,dx\,dy, \int_0^\infty\int_0^\infty \overline{w}(x,y,\alpha)\,dx\,dy\right).$$

***Proof:***

Define the function $I(0, \infty) \to \mathbb{R}^+$ by:

$$\underline{I}(R) = \iint_{D_R} [\underline{w}(x, y, \alpha)] \, dx \, dy, \quad \forall \, \alpha \in [0, 1].$$

From condition 2, it follows that $\underline{I}$ is bounded monotonically increasing. Hence, there exists

$$\lim_{R \to \infty} \underline{I}(R) = \int_0^\infty \int_0^\infty \underline{w}(x, y, \alpha) \, dx \, dy. \quad \blacksquare$$

***Definition (4) [29]:*** Let $\widetilde{w}: (a, b) \times (c, d) \to F^1, x_0 \in (a, b)$. A fuzzy mapping $\widetilde{w}$ is strongly generalized differentiable of the $n^{th}$ order w.r.t. $x$ at $x_0$. If $\widetilde{w}(x, y), \frac{\partial \widetilde{w}(x,y)}{\partial x}, \frac{\partial^2 \widetilde{w}(x,y)}{\partial x^2}, \cdots, \frac{\partial^{(\ell-1)} \widetilde{w}(x,y)}{\partial x^{(\ell-1)}}$ have been strongly generalized differentiable and there exists an element $\frac{\partial^{(\ell)} \widetilde{w}(x_0, y)}{\partial x^{(\ell)}} \in F^1$, $\forall \ell = 1, 2, \cdots, n$:

1. For all $\tau > 0$ is a sufficiently small, there exist $\frac{\partial^{(\ell-1)} \widetilde{w}(x_0 + \tau, y)}{\partial x^{(\ell-1)}} \ominus \frac{\partial^{(\ell-1)} \widetilde{w}(x_0, y)}{\partial x^{(\ell-1)}}, \frac{\partial^{(\ell-1)} \widetilde{w}(x_0, y)}{\partial x^{(\ell-1)}} \ominus \frac{\partial^{(\ell-1)} \widetilde{w}(x_0 - \tau, y)}{\partial x^{(\ell-1)}}$

where

$$\lim_{\tau \to 0} \frac{1}{\tau} \left( \frac{\partial^{(\ell-1)} \widetilde{w}(x_0 + \tau, y)}{\partial x^{(\ell-1)}} \ominus \frac{\partial^{(\ell-1)} \widetilde{w}(x_0, y)}{\partial x^{(\ell-1)}} \right) = \lim_{\tau \to 0} \frac{1}{\tau} \left( \frac{\partial^{(\ell-1)} \widetilde{w}(x_0, y)}{\partial x^{(\ell-1)}} \ominus \frac{\partial^{(\ell-1)} \widetilde{w}(x_0 - \tau, y)}{\partial x^{(\ell-1)}} \right)$$
$$= \frac{\partial^{(\ell)} \widetilde{w}(x_0, y)}{\partial x^{(\ell)}}$$

or

2. For all $\tau > 0$ sufficiently small, there exist $\frac{\partial^{(\ell-1)} \widetilde{w}(x_0, y)}{\partial x^{(\ell-1)}} \ominus \frac{\partial^{(\ell-1)} \widetilde{w}(x_0 + \tau, y)}{\partial x^{(\ell-1)}}$, $\frac{\partial^{(\ell-1)} \widetilde{w}(x_0 - \tau, y)}{\partial x^{(\ell-1)}} \ominus \frac{\partial^{(\ell-1)} \widetilde{w}(x_0, y)}{\partial x^{(\ell-1)}}$

where

$$\lim_{\tau \to 0} \frac{-1}{\tau} \left( \frac{\partial^{(\ell-1)} \widetilde{w}(x_0, y)}{\partial x^{(\ell-1)}} \ominus \frac{\partial^{(\ell-1)} \widetilde{w}(x_0 + \tau, y)}{\partial x^{(\ell-1)}} \right) = \lim_{\tau \to 0} \frac{-1}{\tau} \left( \frac{\partial^{(\ell-1)} \widetilde{w}(x_0 - \tau, y)}{\partial x^{(\ell-1)}} \ominus \frac{\partial^{(\ell-1)} \widetilde{w}(x_0, y)}{\partial x^{(\ell-1)}} \right)$$
$$= \frac{\partial^{(\ell)} \widetilde{w}(x_0, y)}{\partial x^{(\ell)}}$$

or

3. For all $\tau > 0$ is a sufficiently small, there exist $\frac{\partial^{(\ell-1)}\widetilde{w}(x_0+\tau,y)}{\partial x^{(\ell-1)}} \ominus \frac{\partial^{(\ell-1)}\widetilde{w}(x_0,y)}{\partial x^{(\ell-1)}}$, $\frac{\partial^{(\ell-1)}\widetilde{w}(x_0-\tau,y)}{\partial x^{(\ell-1)}} \ominus \frac{\partial^{(\ell-1)}\widetilde{w}(x_0,y)}{\partial x^{(\ell-1)}}$

where

$$\lim_{\tau \to 0} \frac{1}{\tau}\left(\frac{\partial^{(\ell-1)}\widetilde{w}(x_0+\tau,y)}{\partial x^{(\ell-1)}} \ominus \frac{\partial^{(\ell-1)}\widetilde{w}(x_0,y)}{\partial x^{(\ell-1)}}\right) = \lim_{\tau \to 0} \frac{-1}{\tau}\left(\frac{\partial^{(\ell-1)}\widetilde{w}(x_0-\tau,y)}{\partial x^{(\ell-1)}} \ominus \frac{\partial^{(\ell-1)}\widetilde{w}(x_0,y)}{\partial x^{(\ell-1)}}\right)$$
$$= \frac{\partial^{(\ell)}\widetilde{w}(x_0,y)}{\partial x^{(\ell)}}$$

or

4. For all $\tau > 0$ is a sufficiently small, there exist $\frac{\partial^{(\ell-1)}\widetilde{w}(x_0,y)}{\partial x^{(\ell-1)}} \ominus \frac{\partial^{(\ell-1)}\widetilde{w}(x_0+\tau,y)}{\partial x^{(\ell-1)}}$, $\frac{\partial^{(\ell-1)}\widetilde{w}(x_0,y)}{\partial x^{(\ell-1)}} \ominus \frac{\partial^{(\ell-1)}\widetilde{w}(x_0-\tau,y)}{\partial x^{(\ell-1)}}$

where

$$\lim_{\tau \to 0} \frac{-1}{\tau}\left(\frac{\partial^{(\ell-1)}\widetilde{w}(x_0,y)}{\partial x^{(\ell-1)}} \ominus \frac{\partial^{(\ell-1)}\widetilde{w}(x_0+\tau,y)}{\partial x^{(\ell-1)}}\right)$$
$$= \lim_{\tau \to 0} \frac{1}{\tau}\left(\frac{\partial^{(\ell-1)}\widetilde{w}(x_0,y)}{\partial x^{(\ell-1)}} \ominus \frac{\partial^{(\ell-1)}\widetilde{w}(x_0-\tau,y)}{\partial x^{(\ell-1)}}\right) = \frac{\partial^{(\ell)}\widetilde{w}(x_0,y)}{\partial x^{(\ell)}}$$

Similarly,

***Definition (5) [29]:*** Let $\widetilde{w}:(a,b) \times (c,d) \to F^1$, $y_0 \epsilon (a,b)$. A fuzzy mapping $\widetilde{w}$ is strongly generalized differentiable of the $n^{th}$ order w.r.t. $y$ at $y_0$. If $\widetilde{w}(x,y), \frac{\partial \widetilde{w}(x,y)}{\partial y}, \frac{\partial^2 \widetilde{w}(x,y)}{\partial y^2}, \cdots, \frac{\partial^{(\ell-1)}\widetilde{w}(x,y)}{\partial y^{(\ell-1)}}$ have been strongly generalized differentiable and there exists an element $\frac{\partial^{(\ell)}\widetilde{w}(x_0,y)}{\partial y^{(\ell)}} \in F^1$, $\forall \ell = 1, 2, \cdots, n$:

1. For all $\tau > 0$ is a sufficiently small, there exist $\frac{\partial^{(\ell-1)}\widetilde{w}(x,y_0+\tau)}{\partial y^{(\ell-1)}} \ominus \frac{\partial^{(\ell-1)}\widetilde{w}(x,y_0)}{\partial y^{(\ell-1)}}$, $\frac{\partial^{(\ell-1)}\widetilde{w}(x,y_0)}{\partial y^{(\ell-1)}} \ominus \frac{\partial^{(\ell-1)}\widetilde{w}(x,y_0-\tau)}{\partial y^{(\ell-1)}}$

where

$$\lim_{\tau \to 0} \frac{1}{\tau}\left(\frac{\partial^{(\ell-1)}\widetilde{w}(x,y_0+\tau)}{\partial y^{(\ell-1)}} \ominus \frac{\partial^{(\ell-1)}\widetilde{w}(x,y_0)}{\partial y^{(\ell-1)}}\right) = \lim_{\tau \to 0} \frac{1}{\tau}\left(\frac{\partial^{(\ell-1)}\widetilde{w}(x,y_0)}{\partial y^{(\ell-1)}} \ominus \frac{\partial^{(\ell-1)}\widetilde{w}(x,y_0-\tau)}{\partial y^{(\ell-1)}}\right)$$
$$= \frac{\partial^{(\ell)}\widetilde{w}(x,y_0)}{\partial y^{(\ell)}}$$

or

2. For all $\tau > 0$ sufficiently small, there exist $\frac{\partial^{(\ell-1)}\tilde{w}(x,y_0)}{\partial y^{(\ell-1)}} \ominus \frac{\partial^{(\ell-1)}\tilde{w}(x,y_0+\tau)}{\partial y^{(\ell-1)}}$, $\frac{\partial^{(\ell-1)}\tilde{w}(x,y_0-\tau)}{\partial y^{(\ell-1)}} \ominus \frac{\partial^{(\ell-1)}\tilde{w}(x,y_0)}{\partial y^{(\ell-1)}}$

where

$$\lim_{\tau \to 0} \frac{-1}{\tau} \left( \frac{\partial^{(\ell-1)}\tilde{w}(x, y_0)}{\partial y^{(\ell-1)}} \ominus \frac{\partial^{(\ell-1)}\tilde{w}(x, y_0 + \tau)}{\partial y^{(\ell-1)}} \right) = \lim_{\tau \to 0} \frac{-1}{\tau} \left( \frac{\partial^{(\ell-1)}\tilde{w}(x, y_0 - \tau)}{\partial y^{(\ell-1)}} \ominus \frac{\partial^{(\ell-1)}\tilde{w}(x, y_0)}{\partial y^{(\ell-1)}} \right)$$
$$= \frac{\partial^{(\ell)}\tilde{w}(x, y_0)}{\partial y^{(\ell)}}$$

or

3. For all $\tau > 0$ is a sufficiently small, there exist $\frac{\partial^{(\ell-1)}\tilde{w}(x,y_0+\tau)}{\partial y^{(\ell-1)}} \ominus \frac{\partial^{(\ell-1)}\tilde{w}(x,y_0)}{\partial y^{(\ell-1)}}$, $\frac{\partial^{(\ell-1)}\tilde{w}(x,y_0-\tau)}{\partial y^{(\ell-1)}} \ominus \frac{\partial^{(\ell-1)}\tilde{w}(x,y_0)}{\partial y^{(\ell-1)}}$

where

$$\lim_{\tau \to 0} \frac{1}{\tau} \left( \frac{\partial^{(\ell-1)}\tilde{w}(x, y_0 + \tau)}{\partial y^{(\ell-1)}} \ominus \frac{\partial^{(\ell-1)}\tilde{w}(x, y_0)}{\partial y^{(\ell-1)}} \right) = \lim_{\tau \to 0} \frac{-1}{\tau} \left( \frac{\partial^{(\ell-1)}\tilde{w}(x, y_0 - \tau)}{\partial y^{(\ell-1)}} \ominus \frac{\partial^{(\ell-1)}\tilde{w}(x, y_0)}{\partial y^{(\ell-1)}} \right)$$
$$= \frac{\partial^{(\ell)}\tilde{w}(x, y_0)}{\partial y^{(\ell)}}$$

or

4. For all $\tau > 0$ is a sufficiently small, there exist $\frac{\partial^{(\ell-1)}\tilde{w}(x,y_0)}{\partial y^{(\ell-1)}} \ominus \frac{\partial^{(\ell-1)}\tilde{w}(x,y_0+\tau)}{\partial y^{(\ell-1)}}$, $\frac{\partial^{(\ell-1)}\tilde{w}(x,y_0)}{\partial y^{(\ell-1)}} \ominus \frac{\partial^{(\ell-1)}\tilde{w}(x,y_0-\tau)}{\partial y^{(\ell-1)}}$

where

$$\lim_{\tau \to 0} \frac{-1}{\tau} \left( \frac{\partial^{(\ell-1)}\tilde{w}(x, y_0)}{\partial y^{(\ell-1)}} \ominus \frac{\partial^{(\ell-1)}\tilde{w}(x, y_0 + \tau)}{\partial y^{(\ell-1)}} \right) = \lim_{\tau \to 0} \frac{1}{\tau} \left( \frac{\partial^{(\ell-1)}\tilde{w}(x, y_0)}{\partial y^{(\ell-1)}} \ominus \frac{\partial^{(\ell-1)}\tilde{w}(x, y_0 - \tau)}{\partial y^{(\ell-1)}} \right)$$
$$= \frac{\partial^{(\ell)}\tilde{w}(x, y_0)}{\partial y^{(\ell)}}$$

***Theorem (2) [29]:*** Let $\tilde{w}(x,y), \frac{\partial \tilde{w}(x,y)}{\partial x}, \frac{\partial^2 \tilde{w}(x,y)}{\partial x^2}, \frac{\partial^3 \tilde{w}(x,y)}{\partial x^3}, \cdots, \frac{\partial^{(n-1)}\tilde{w}(x,y)}{\partial x^{(n-1)}}$ are differentiable fuzzy-valued functions on $[0, \infty)$. A fuzzy function which represented $\tilde{w}(x, y) = [\underline{w}(x, y, \alpha), \overline{w}(x, y, \alpha)]$ for all $\alpha \in [0, 1]$, then:

$$\frac{\partial^{(n)}\widetilde{w}(x,y)}{\partial x^{(n)}} = \begin{cases} \left[\dfrac{\partial^{(n)}\underline{w}(x,y,\alpha)}{\partial x^{(n)}}, \dfrac{\partial^{(n)}\overline{w}(x,y,\alpha)}{\partial x^{(n)}}\right] & \text{if number of (ii) − differentiable is even,} \\ \left[\dfrac{\partial^{(n)}\overline{w}(x,y,\alpha)}{\partial x^{(n)}}, \dfrac{\partial^{(n)}\underline{w}(x,y,\alpha)}{\partial x^{(n)}}\right] & \text{if number of (ii) − differentiable is odd.} \end{cases}$$

**_Theorem (3) [29]:_** Let $\widetilde{w}(x,y), \frac{\partial \widetilde{w}(x,y)}{\partial y}, \frac{\partial^2 \widetilde{w}(x,y)}{\partial y^2}, \frac{\partial^3 \widetilde{w}(x,y)}{\partial y^3}, \cdots, \frac{\partial^{(n-1)} \widetilde{w}(x,y)}{\partial y^{(n-1)}}$ are differentiable fuzzy-valued functions on $[0, \infty)$. A fuzzy function which represented $\widetilde{w}(x,y) = [\underline{w}(x,y,\alpha), \overline{w}(x,y,\alpha)]$ for all $\alpha \in [0, 1]$, then:

$$\frac{\partial^{(n)}\widetilde{w}(x,y)}{\partial y^{(n)}} = \begin{cases} \left[\dfrac{\partial^{(n)}\underline{w}(x,y,\alpha)}{\partial y^{(n)}}, \dfrac{\partial^{(n)}\overline{w}(x,y,\alpha)}{\partial y^{(n)}}\right] & \text{if number of (ii) − differentiable is even,} \\ \left[\dfrac{\partial^{(n)}\overline{w}(x,y,\alpha)}{\partial y^{(n)}}, \dfrac{\partial^{(n)}\underline{w}(x,y,\alpha)}{\partial y^{(n)}}\right] & \text{if number of (ii) − differentiable is odd.} \end{cases}$$

**_Theorem (4) [27]:_** let $\widetilde{w}: D \to F^1$ be fuzzy-valued function with parametric form $\left(\underline{w}(x,y,\alpha), \overline{w}(x,y,\alpha)\right)$ for all $\alpha \in [0, 1]$:

1. The functions $\underline{w}(x,y,\alpha)$ and $\overline{w}(x,y,\alpha)$ are Riemann-integrable on $D_\mathbb{R}$.
2. There are constants $\underline{M}(\alpha) > 0$ and $\overline{M}(\alpha) > 0$, such that:

$$\iint_{D_R} [\underline{w}(x,y,\alpha)]\, dx\, dy \leq \underline{M}(\alpha)$$

$$\iint_{D_R} [\overline{w}(x,y,\alpha)]\, dx\, dy \leq \overline{M}(\alpha), \quad \forall\, R > 0,$$

Then the function $\widetilde{w}(x,y)$ is an improper fuzzy Riemann-integrable on $D$ and

$$(FR)\int_0^\infty (FR)\int_0^\infty \widetilde{w}(x,y)\, dx\, dy = \left(\int_0^\infty \int_0^\infty \underline{w}(x,y,\alpha)\, dx\, dy, \int_0^\infty \int_0^\infty \overline{w}(x,y,\alpha)\, dx\, dy\right).$$

**_Proof:_**

Define the function $I(0, \infty) \to \mathbb{R}^+$ by:

$$\underline{I}(R) = \iint_{D_R} [\underline{w}(x,y,\alpha)]\, dx\, dy, \quad \forall\, \alpha \in [0, 1].$$

From condition 2, it follows that $\underline{I}$ is bounded monotonically increasing. Hence, there exists

$$\lim_{R\to\infty} \underline{I}(R) = \int_0^\infty \int_0^\infty \underline{w}(x,y,\alpha)\,dx\,dy. \quad \blacksquare$$

**_Proposition (1) [29]:_** If $\tilde{w}(x,y)$ is defined on fuzzy Riemann-integrable on $D_R$ and if $\tilde{w}(x,y) = \left(\underline{w}(x,y,\alpha), \overline{w}(x,y,\alpha)\right)$ such that $\underline{w}(x,y,\alpha) \leq \overline{w}(x,y,\alpha)$ for all $x,y \in D_R$, and $\alpha \in [0,1]$ we have

$$\iint_{D_R} \underline{w}(x,y,\alpha)\,dx\,dy \leq \iint_{D_R} \overline{w}(x,y,\alpha)\,dx\,dy.$$

In 2022 Ahmad Issa, Emad A. Kuffi presented an integral transform named Complex EE transform. This transform is defined for the function f(t) as: [28]

$$D_{EE}[f(x,t)] = A(iu, iv) = \int_0^\infty \int_0^\infty f(x,t) e^{-i(u^n x + v^n t)}\,dt\,dx, n \in Z, t, x \geq 0, i \in C.$$

So the definition of DFCEET will depend on the above transform as follows:

**_Definition (6):_** let $\tilde{w}: D \to F^1$ be a continuous fuzzy-valued function and the function $e^{-i(u^n x + v^n y)} \odot \tilde{w}(x,y)$ is improper fuzzy Riemann-integrable on $D$, then

$$(F\mathcal{R})\int_0^\infty (F\mathcal{R})\int_0^\infty e^{-i(u^n x + v^n y)} \odot \tilde{w}(x,y)\,dx\,dy. \quad (1)$$

where $n \in Z, y, x \geq 0, i \in C.$

So we can define fuzzy double EE integral transform:

$$\tilde{D}_{EE}[w(x,y)] = \tilde{A}(iu, iv) = \int_0^\infty \int_0^\infty e^{-i(u^n x + v^n y)} \odot \tilde{w}(x,y)\,dx\,dy. \quad (2)$$

The parametric form fuzzy EE double integral transform is follows:

$$\tilde{D}_{EE}[w(x,y)] = \left(\underline{D}_{EE}[w(x,y)], \overline{D}_{EE}[w(x,y)]\right) = \left(\underline{A}(iu, iv, \alpha), \overline{A}(iu, iv, \alpha)\right), \quad (3)$$

where

$$\underline{D}_{EE}[w(x,y)] = \underline{A}(iu, iv, \alpha) = \int_0^\infty \int_0^\infty e^{-i(u^n x + v^n y)} \underline{w}(x, y, \alpha) \, dx \, dy,$$

$$\overline{D}_{EE}[w(x,y)] = \overline{A}(iu, iv, \alpha) = \int_0^\infty \int_0^\infty e^{-i(u^n x + v^n y)} \overline{w}(x, y, \alpha) \, dx \, dy.$$

***Definition (7):*** The inverse of fuzzy EE double integral transform can be written as the formula

$$\widetilde{D}_{EE}^{-1}\{A(iu, iv)\} = \frac{1}{2\pi i} \int_{\gamma-i\infty}^{\gamma+i\infty} ie^{iu^n x} ds \left( \frac{1}{2\pi i} \int_{\delta-i\infty}^{\delta+i\infty} ie^{iv^n y} \widetilde{A}(iu, iv) dr \right),$$

it follows

$$\widetilde{D}_{EE}^{-1}\{A(iu, iv)\} = \left( \underline{D}_{EE}^{-1}[A(iu, iv)], \overline{D}_{EE}^{-1}[A(iu, iv)] \right),$$

Where

$$\underline{D}_{EE}^{-1}[A(iu, iv)] = \frac{1}{2\pi i} \int_{\gamma-i\infty}^{\gamma+i\infty} ie^{iu^n x} ds \left( \frac{1}{2\pi i} \int_{\delta-i\infty}^{\delta+i\infty} ie^{iv^n y} \underline{A}(iu, iv) dr \right),$$

and

$$\overline{D}_{EE}^{-1}[A(iu, iv)] = \frac{1}{2\pi i} \int_{\gamma-i\infty}^{\gamma+i\infty} ie^{iu^n x} ds \left( \frac{1}{2\pi i} \int_{\delta-i\infty}^{\delta+i\infty} ie^{iv^n y} \overline{A}(iu, iv) dr \right),$$

for all $\alpha \in [0, 1]$. The functions $\underline{A}(iu, iv, \alpha)$ and $\overline{A}(iu, iv, \alpha)$ must be analytic functions for all $u^n$ and $v^n$ in the region defined by the inequalities $\Re e(iu^n) \geq \gamma$ and $\Re e(iv^n) \geq \delta$, where $\gamma, \delta$ are real constants to be chosen suitably.

The Double Complex EE Transform is applied for some elementary functions [28]:

1. $D_{EE}\{1\} = \frac{-i}{u^n} \cdot \frac{-i}{v^n}$.
2. $D_{EE}\{x^r t^m\} = \frac{r!}{(iu)^{(r+1)n}} \cdot \frac{m!}{(iv)^{(m+1)n}}$.
3. $D_{EE}\{e^{ax+bt}\} = \frac{(a+iu^n)(b+iv^n)}{(a^2+u^{2n})(b^2+v^{2n})}$.
4. $D_{EE}\{e^{-(ax+bt)}\} = \frac{(a-iu^n)(b-iv^n)}{(a^2+u^{2n})(b^2+v^{2n})}$.
5. $D_{EE}\{e^{i(ax+bt)}\} = \frac{-1}{(u^n-a)(v^n-b)}$.
6. $D_{EE}\{e^{-i(ax+bt)}\} = \frac{-1}{(u^n+a)(v^n+b)}$.
7. $D_{EE}\{\sin(ax+bt)\} = \frac{i(bu^n+av^n)}{(u^{2n}-a^2)(v^{2n}-b^2)}$.
8. $D_{EE}\{\cos(ax+bt)\} = \frac{-((uv)^n+ab)}{(u^{2n}-a^2)(v^{2n}-b^2)}$.
9. $D_{EE}\{\sinh(ax+bt)\} = \frac{i(av^n+bu^n)}{(a^2+u^{2n})(b^2+v^{2n})}$.

10. $D_{EE}\{cosh(ax+bt)\} = \frac{ab-(uv)^n}{(a^2+u^{2n})(b^2+v^{2n})}$.

**_Theorem (5):_** let $\tilde{w}(x,y)$ be a continuous fuzzy-valued function. If $\tilde{D}_{EE}\{w(x,y)\}$ is a fuzzy double EE integral transform of $\tilde{w}(x,y)$, then

$$\tilde{D}_{EE}\{e^{ax+by} \odot w(x,y)\} = \tilde{A}\big((iu^n - a), (iv^n - b)\big),$$

where $a, b$ are arbitrary constants.

**_Proof:_** we have

$$\tilde{D}_{EE}\{e^{ax+by} \odot w(x,y)\} = \big(\underline{D}_{EE}\{e^{ax+by}\underline{w}(x,y,\alpha)\}, \overline{D}_{EE}\{e^{ax+by}\overline{w}(x,y,\alpha)\}\big),$$

$$= \left( \int_0^\infty \int_0^\infty e^{ax+by} e^{-i(u^n x + v^n y)} \underline{w}(x,y,\alpha)\, dx\, dy, \right.$$

$$\left. \int_0^\infty \int_0^\infty e^{ax+by} e^{-i(u^n x + v^n y)} \overline{w}(x,y,\alpha)\, dx\, dy. \right)$$

$$\tilde{D}_{EE}\{e^{ax+by} \odot w(x,y)\}$$

$$= \left( \int_0^\infty \int_0^\infty e^{-((iu^n-a)x+(iv^n-b)y)} \underline{w}(x,y,\alpha)\, dx\, dy, \right.$$

$$\left. \int_0^\infty \int_0^\infty e^{-((iu^n-a)x+(iv^n-b)y)} \overline{w}(x,y,\alpha)\, dx\, dy \right),$$

$$= \big(\underline{A}((iu^n-a),(iv^n-b),\alpha), \overline{A}((iu^n-a),(iv^n-b),\alpha)\big),$$

$$= \tilde{A}\big((iq_1(s) - a), (iq_2(r) - b)\big). \quad \blacksquare$$

**_Definition (8) [27]:_** If $\tilde{K}: D \to F^1$ and $\tilde{w}(x,y)$ are fuzzy Riemann-integrable functions, then the double fuzzy convolution $\tilde{K}(x,y)$ and $\tilde{w}(x,y)$ is given by:

$$(\tilde{K} ** \tilde{w})(x,y) = (FR)\int_0^y (FR)\int_0^x \tilde{K}(x-\tau, y-\mu)\tilde{w}(\tau,\mu)\, d\tau\, d\mu.$$

The symbol $**$ denotes the double convolution respect to $x$ and $y$.

**_Theorem (6):_** let $\tilde{K}: D \to F^1$ and $\tilde{w}(x,y)$ be continuous fuzzy-valued functions. Then the fuzzy double EE integral transform of the double fuzzy convolution $\tilde{K}$ and $\tilde{w}$ is given by:

$$\widetilde{D}_{EE}\{(\widetilde{K} ** \widetilde{w})(x,y)\} = \widetilde{D}_{EE}\{K(x,y)\} \odot \widetilde{D}_{EE}\{w(x,y)\} \quad (4)$$

***Proof:*** we have

$$\widetilde{D}_{EE}\{(\widetilde{K} ** \widetilde{w})(x,y)\} =$$

$$\int_0^\infty \int_0^\infty e^{-i(u^n x + v^n y)} \odot \left\{ \int_0^y \int_0^x \widetilde{K}(x-\tau, y-\mu) \widetilde{w}(\tau,\mu) \, d\tau \, d\mu \right\} dx \, dy,$$

$$= \int_0^\infty \int_0^\infty \widetilde{w}(\tau,\mu) \, d\tau \, d\mu \int_0^\infty \int_0^\infty e^{-i(u^n x + v^n y)} \odot \widetilde{K}(x-\tau, y-\mu) \, dx \, dy,$$

Now, we assume that $x - \tau = m$ and $-\mu = l$, we have:

$$\widetilde{D}_{EE}\{(\widetilde{K} ** \widetilde{w})(x,y)\} = \int_0^\infty \int_0^\infty \widetilde{w}(\tau,\mu) \, d\tau \, d\mu \int_0^\infty \int_0^\infty e^{-i[u^n(m+\tau) + v^n(l+\mu)]} \odot \widetilde{K}(m,l) \, dm \, dl,$$

$$= \int_0^\infty \int_0^\infty e^{-i(u^n \tau + v^n \mu)} \odot \widetilde{w}(\tau,\mu) \, d\tau \, d\mu \int_0^\infty \int_0^\infty e^{-i(u^n m + v^n l)} \odot \widetilde{K}(m,l) \, dm \, dl,$$

$$\widetilde{D}_{EE}\{(\widetilde{K} ** \widetilde{w})(x,y)\} = \widetilde{D}_{EE}\{K(x,y)\} \odot \widetilde{D}_{EE}\{w(x,y)\}. \quad \blacksquare$$

***Theorem (7):*** let $\widetilde{w}: D \to F^1$ be a continuous fuzzy-valued function and the function $e^{-i(u^n x + v^n y)} \odot \widetilde{w}(x,y)$ and $e^{-i(u^n x + v^n y)} \odot \frac{\partial^j \widetilde{w}(x,y)}{\partial x^j}$ are improper fuzzy Riemann-integrable on $D$, then

$$\widetilde{D}_{EE}\left\{\frac{\partial^j \widetilde{w}(x,y)}{\partial x^j}\right\} = \frac{\partial^j}{\partial x^j} \widetilde{D}_{EE}\{\widetilde{w}(x,y)\}.$$

***Proof:***

Let $\widetilde{w}(x,y)$ be (i)-differentiable. From (1)

$$\widetilde{D}_{EE}\left\{\frac{\partial^j \widetilde{w}(x,y)}{\partial x^j}\right\} = (F\mathcal{R}) \int_0^\infty (F\mathcal{R}) \int_0^\infty e^{-i(u^n x + v^n y)} \odot \frac{\partial^j \widetilde{w}(x,y)}{\partial x^j} \, dx \, dy,$$

$$= \left( \int_0^\infty \int_0^\infty e^{-i(u^n x + v^n y)} \frac{\partial^j \underline{w}(x,y)}{\partial x^j} \, dx \, dy, \int_0^\infty \int_0^\infty e^{-i(u^n x + v^n y)} \frac{\partial^j \overline{w}(x,y)}{\partial x^j} \, dx \, dy \right),$$

$$= \frac{\partial^j}{\partial x^j}\left(\int_0^\infty\int_0^\infty e^{-i(u^n x+v^n y)} \underline{w}(x,y)dx\,dy, \int_0^\infty\int_0^\infty e^{-i(u^n x+v^n y)}\overline{w}(x,y)\,dx\,dy\right),$$

$$= \frac{\partial^j}{\partial x^j}\widetilde{D}_{EE}\{\widetilde{w}(x,y)\}. \quad \blacksquare$$

***Theorem (8):*** let $\widetilde{w}: D \to F^1$ be a continuous fuzzy-valued function and the function $e^{-i(u^n x+v^n y)} \odot \widetilde{w}(x,y)$ and $e^{-i(u^n x+v^n y)} \odot \frac{\partial^j \widetilde{w}(x,y)}{\partial y^j}$ are improper fuzzy Riemann-integrable on $D$, then

$$\widetilde{D}_{EE}\left\{\frac{\partial^j \widetilde{w}(x,y)}{\partial y^j}\right\} = \frac{\partial^j}{\partial y^j}\widetilde{D}_{EE}\{\widetilde{w}(x,y)\}.$$

***Proof:***

Let $\widetilde{w}(x,y)$ be (i)-differentiable. From (1)

$$\widetilde{D}_{EE}\left\{\frac{\partial^j \widetilde{w}(x,y)}{\partial y^j}\right\} = (F\mathcal{R})\int_0^\infty (F\mathcal{R})\int_0^\infty e^{-i(u^n x+v^n y)} \odot \frac{\partial^j \widetilde{w}(x,y)}{\partial y^j}\,dx\,dy,$$

$$= \left(\int_0^\infty\int_0^\infty e^{-i(u^n x+v^n y)} \frac{\partial^j \underline{w}(x,y)}{\partial y^j}dx\,dy, \int_0^\infty\int_0^\infty e^{-i(u^n x+v^n y)} \frac{\partial^j \overline{w}(x,y)}{\partial y^j}dx\,dy\right),$$

$$= \frac{\partial^j}{\partial y^j}\left(\int_0^\infty\int_0^\infty e^{-i(u^n x+v^n y)} \underline{w}(x,y)dx\,dy, \int_0^\infty\int_0^\infty e^{-i(u^n x+v^n y)}\overline{w}(x,y)\,dx\,dy\right),$$

$$= \frac{\partial^j}{\partial y^j}\widetilde{D}_{EE}\{\widetilde{w}(x,y)\}. \quad \blacksquare$$

***Theorem (9):*** let $\widetilde{w}: D \to F^1$ be a continuous fuzzy-valued function and the function $e^{-i(u^n x+v^n y)} \odot \widetilde{w}(x,y)$ and $e^{-i(u^n x+v^n y)} \odot \frac{\partial^j \widetilde{w}(x,y)}{\partial x^j}$ are improper fuzzy Riemann-integrable on $D$. For all $x > 0$ and $j \in \mathbb{N}$ there exist to continuous partial H-derivatives to $(i-1)-th$ order with respect to $x$ and there exists $\frac{\partial^j \widetilde{w}(x,y)}{\partial x^j}$. Then

$$\widetilde{D}_{EE}\left\{\frac{\partial^j \widetilde{w}(x,y)}{\partial x^j}\right\} = (iu^n)^j \odot \widetilde{D}_{EE}\{w(x,y)\} - \left[\sum_{k=0}^{j-1}(iu^n)^{j-k-1} \odot \widetilde{E}^c\left\{\frac{\partial^k w}{\partial x^k}(0,y)\right\}\right]$$

1. If the function $\widetilde{w}(x,y)$ is (i)-differentiable then
$$\widetilde{D}_{EE}\left\{\frac{\partial^j \widetilde{w}(x,y)}{\partial x^j}\right\} = \left(\underline{D}_{EE}\left[\frac{\partial^j \underline{w}(x,y,\alpha)}{\partial x^j}\right], \overline{D}_{EE}\left[\frac{\partial^j \overline{w}(x,y,\alpha)}{\partial x^j}\right]\right),$$

2. If the function $\widetilde{w}(x,y)$ is (ii)-differentiable then
$$\widetilde{D}_{EE}\left\{\frac{\partial^j \widetilde{w}(x,y)}{\partial x^j}\right\} = \left(\overline{D}_{EE}\left[\frac{\partial^j \underline{w}(x,y,\alpha)}{\partial x^j}\right], \underline{D}_{EE}\left[\frac{\partial^j \overline{w}(x,y,\alpha)}{\partial x^j}\right]\right),$$

Where

$$\underline{D}_{EE}\left[\frac{\partial^j \underline{w}(x,y,\alpha)}{\partial x^j}\right] = (iu^n)^j \underline{D}_{EE}\{\underline{w}(x,y,\alpha)\} - \left[\sum_{k=0}^{j-1}(iu^n)^{j-k-1}\underline{E}^c\left\{\frac{\partial^k \underline{w}}{\partial x^k}(0,y,\alpha)\right\}\right],$$

$$\overline{D}_{EE}\left[\frac{\partial^j \overline{w}(x,y,\alpha)}{\partial x^j}\right]$$
$$= (iu^n)^j \overline{D}_{EE}\{\overline{w}(x,y,\alpha)\} - \left[\sum_{k=0}^{j-1}(iu^n)^{j-k-1}\overline{E^c}\left\{\frac{\partial^k \overline{w}}{\partial x^k}(0,y,\alpha)\right\}\right]. \quad (5)$$

***Proof:***

Let the function $\widetilde{w}(x,y)$ is (i)-differentiable. By induction, for $n=1$, from the condition (3) we have

$$\widetilde{D}_{EE}\left\{\frac{\partial \widetilde{w}(x,y)}{\partial x}\right\} = \left(\underline{D}_{EE}\left\{\frac{\partial \underline{w}(x,y,\alpha)}{\partial x}\right\}, \overline{D}_{EE}\left\{\frac{\partial \overline{w}(x,y,\alpha)}{\partial x}\right\}\right).$$

By us part integration on $x$ and condition (3) we obtain

$$\widetilde{D}_{EE}\left\{\frac{\partial \widetilde{w}(x,y)}{\partial x}\right\} = \left(\int_0^\infty\int_0^\infty e^{-i(u^n x+v^n y)}\frac{\partial \underline{w}(x,y)}{\partial x}dx\,dy, \int_0^\infty\int_0^\infty e^{-i(u^n x+v^n y)}\frac{\partial \overline{w}(x,y)}{\partial x}dx\,dy\right),$$

$$= \left(iu^n \underline{D}_{EE}\left\{\frac{\partial \underline{w}(x,y,\alpha)}{\partial x}\right\} - \underline{E}^c\left\{\frac{\partial \underline{w}(0,y,\alpha)}{\partial x}\right\}, iu^n \overline{D}_{EE}\left\{\frac{\partial \overline{w}(x,y,\alpha)}{\partial x}\right\}\right.$$
$$\left. - \overline{E^c}\left\{\frac{\partial \overline{w}(0,y,\alpha)}{\partial x}\right\}\right).$$

For $j=m$, the equations (5) are hold, then

$$\widetilde{D}_{EE}\left\{\frac{\partial^m \widetilde{w}(x,y)}{\partial x^m}\right\} = \left((iu^n)^m \underline{D}_{EE}\{\underline{w}(x,y,\alpha)\} - \left[\sum_{k=0}^{m-1}(iu^n)^{m-k-1}\underline{E}^c\left\{\frac{\partial^k \underline{w}}{\partial x^k}(0,y,\alpha)\right\}\right],\right.$$

$$(iu^n)^m \overline{D}_{EE}\{w(x,y,\alpha)\} - \left[\sum_{k=0}^{m-1}(iu^n)^{m-k-1}\overline{E^c}\left\{\frac{\partial^k w}{\partial x^k}(0,y,\alpha)\right\}\right]).$$

Hence, for $j = m + 1$ we get

$$\widetilde{D}_{EE}\left\{\frac{\partial^{m+1}w(x,y)}{\partial x^{m+1}}\right\} = \left(\underline{D}_{EE}\left\{\frac{\partial^{m+1}w}{\partial x^{m+1}}(x,y,\alpha)\right\}, \overline{D}_{EE}\left\{\frac{\partial^{m+1}w}{\partial x^{m+1}}(x,y,\alpha)\right\}\right),$$

$$\frac{\partial}{\partial x}\widetilde{D}_{EE}\left\{\frac{\partial^m}{\partial x^m}w(x,y)\right\} = \left(\frac{\partial}{\partial x}\underline{D}_{EE}\left\{\frac{\partial^m w}{\partial x^m}(x,y,\alpha)\right\}, \frac{\partial}{\partial x}\overline{D}_{EE}\left\{\frac{\partial^m w}{\partial x^m}(x,y,\alpha)\right\}\right),$$

$$= \left(\frac{\partial}{\partial x}\left((iu^n)^m \underline{D}_{EE}\{w(x,y,\alpha)\} - \left[\sum_{k=0}^{m-1}(iu^n)^{m-k-1}\underline{E^c}\left\{\frac{\partial^k w}{\partial x^k}(0,y,\alpha)\right\}\right]\right),$$

$$\frac{\partial}{\partial x}\left((iu^n)^m \overline{D}_{EE}\{w(x,y,\alpha)\} - \left[\sum_{k=0}^{m-1}(iu^n)^{m-k-1}\overline{E^c}\left\{\frac{\partial^k w}{\partial x^k}(0,y,\alpha)\right\}\right]\right)\right),$$

$$= \left((iu^n)^{m+1}\left[\underline{D}_{EE}\{w(x,y,\alpha)\} - \underline{E^c}\{w(0,y,\alpha)\}\right] - \left[\sum_{k=0}^{m-1}(iu^n)^{m-k-1}\underline{E^c}\left\{\frac{\partial^k w}{\partial x^k}(0,y,\alpha)\right\}\right],$$

$$(iu^n)^{m+1}\left[\overline{D}_{EE}\{w(x,y,\alpha)\} - \overline{E^c}\{w(0,y,\alpha)\}\right] - \left[\sum_{k=0}^{m-1}(iu^n)^{m-k-1}\overline{E^c}\left\{\frac{\partial^k w}{\partial x^k}(0,y,\alpha)\right\}\right]\right)$$

$$= \left((iu^n)^{m+1}\underline{D}_{EE}\{w(x,y,\alpha)\} - \left[\sum_{k=0}^{m}(iu^n)^{m-k-1}\underline{E^c}\left\{\frac{\partial^k w}{\partial x^k}(0,y,\alpha)\right\}\right],$$

$$(iu^n)^{m+1}\overline{D}_{EE}\{w(x,y,\alpha)\} - \left[\sum_{k=0}^{m}(iu^n)^{m-k-1}\overline{E^c}\left\{\frac{\partial^k w}{\partial x^k}(0,y,\alpha)\right\}\right]\right). \blacksquare$$

Similarly, we can prove the following theorem

**_Theorem (10):_** let $\widetilde{w}: D \to F^1$ be a continuous fuzzy-valued function and the function $e^{-i(u^n x + v^n y)} \odot \widetilde{w}(x,y)$ and $e^{-i(u^n x + v^n y)} \odot \frac{\partial^j \widetilde{w}(x,y)}{\partial y^j}$ are improper fuzzy Riemann-integrable on $D$. For all $y > 0$ and $j \in \mathbb{N}$ there exist to continuous partial H-derivatives to $(i-1)-th$ order with respect to $y$ and there exists $\frac{\partial^j \widetilde{w}(x,y)}{\partial y^j}$. Then

$$\widetilde{D}_{EE}\left\{\frac{\partial^j \widetilde{w}(x,y)}{\partial y^j}\right\} = (iv^n)^j \odot \widetilde{D}_{EE}\{w(x,y)\} - \left[\sum_{k=0}^{j-1}(iv^n)^{j-k-1}\odot \widetilde{E}^c\left\{\frac{\partial^k w}{\partial y^k}(x,0)\right\}\right]$$

1. If the function $\widetilde{w}(x,y)$ is (i)-differentiable then

$$\widetilde{D}_{EE}\left\{\frac{\partial^j \widetilde{w}(x,y)}{\partial y^j}\right\} = \left(\underline{D}_{EE}\left[\frac{\partial^j w(x,y,\alpha)}{\partial y^j}\right], \overline{D}_{EE}\left[\frac{\partial^j w(x,y,\alpha)}{\partial y^j}\right]\right),$$

2. If the function $\widetilde{w}(x,y)$ is (ii)-differentiable then

$$\widetilde{D}_{EE}\left\{\frac{\partial^j \widetilde{w}(x,y)}{\partial y^j}\right\} = \left(\overline{D}_{EE}\left[\frac{\partial^j w(x,y,\alpha)}{\partial y^j}\right], \underline{D}_{EE}\left[\frac{\partial^j w(x,y,\alpha)}{\partial y^j}\right]\right),$$

Where

$$\underline{D}_{EE}\left\{\frac{\partial^j w(x,y,\alpha)}{\partial y^j}\right\} = (iv^n)^j \underline{D}_{EE}\{w(x,y,\alpha)\} - \left[\sum_{k=0}^{j-1}(iv^n)^{j-k-1}\underline{E^c}\left\{\frac{\partial^k w}{\partial y^k}(x,0,\alpha)\right\}\right], \quad (6)$$

$$\overline{D}_{EE}\left\{\frac{\partial^j w(x,y,\alpha)}{\partial y^j}\right\} = (iv^n)^j \overline{D}_{EE}\{w(x,y,\alpha)\} - \left[\sum_{k=0}^{j-1}(iv^n)^{j-k-1}\overline{E^c}\left\{\frac{\partial^k w}{\partial y^k}(x,0,\alpha)\right\}\right]. \quad (7)$$

***Double Fuzzy Complex EE Transform for Solving Partial Voltera Fuzzy Integro-Differential Equations***

In this section, we apply the double fuzzy complex EE transform to solve partial voltera fuzzy integro-differential equations. This equation is defined as

$$\sum_{h=1}^{l} a_h \odot \frac{\partial^h \widetilde{w}(x,y)}{\partial x^h} \oplus \sum_{j=1}^{m} b_j \odot \frac{\partial^j \widetilde{w}(x,y)}{\partial y^j} \oplus c \odot \widetilde{w}(x,y)$$
$$= \widetilde{g}(x,y) \oplus (FR)\int_0^y (FR)\int_0^x K(x-\tau, y-\mu) \odot \widetilde{w}(\tau,\mu)\, d\tau\, d\mu, \quad (8)$$

With initial conditions

$$\frac{\partial^h \widetilde{w}(0,y)}{\partial x^h} = \Phi_h(y), \quad h = 0,1,\ldots,l-1,$$

$$\frac{\partial^j \widetilde{w}(0,y)}{\partial y^j} = \Psi_j(x), \quad j = 0,1,\ldots,m-1,$$

Where $K:[0,b] \times [0,d] \to \mathbb{R}$, is a continuous function, and $\widetilde{g}, \widetilde{w} : [0,b] \times [0,d] \to F^1, \Phi_h : [0,d] \to F^1, \Psi_j : [0,b] \to F^1$ are continuous fuzzy functions and $a_h, h = 0,1,\ldots,l, b_j, j = 0,1,\ldots,m, c$, are constants.

Applying (2) on both sides of (8) to get

$$\widetilde{D}_{EE}\left\{\sum_{h=1}^{l} a_h \odot \frac{\partial^h \widetilde{w}(x,y)}{\partial x^h}\right\} \oplus \widetilde{D}_{EE}\left\{\sum_{j=1}^{m} b_j \odot \frac{\partial^j \widetilde{w}(x,y)}{\partial y^j}\right\} \oplus \widetilde{D}_{EE}\{c \odot \widetilde{w}(x,y)\}$$
$$= \widetilde{D}_{EE}\{\widetilde{g}(x,y)\} \oplus \widetilde{D}_{EE}\left\{(FR)\int_0^y (FR)\int_0^x K(x-\tau, y-\mu) \odot \widetilde{w}(\tau,\mu)\, d\tau\, d\mu\right\},$$

By using fuzzy convolution (4) we obtain

$$\sum_{h=1}^{l} a_h \odot \widetilde{D}_{EE}\left\{\frac{\partial^h w(x,y)}{\partial x^h}\right\} \oplus \sum_{j=1}^{m} b_j \odot \widetilde{D}_{EE}\left\{\frac{\partial^j w(x,y)}{\partial y^j}\right\} \oplus c \odot \widetilde{D}_{EE}\{w(x,y)\}$$
$$= \widetilde{D}_{EE}\{g(x,y)\} \oplus D_{EE}\{K(x,y)\} \odot \widetilde{D}_{EE}\{w(x,y)\},$$

Let the constants $a_h, h = 0, 1, \ldots, l, b_j, j = 0, 1, \ldots, m, c$ be positive and the function $K(x,y) > 0$.

1. If $\widetilde{w}(x,y)$ be (i)-differentiable, then

$$\sum_{h=1}^{l} a_h \underline{D}_{EE}\left\{\frac{\partial^h w(x,y)}{\partial x^h}\right\} + \sum_{j=1}^{m} b_j \underline{D}_{EE}\left\{\frac{\partial^j w(x,y)}{\partial y^j}\right\} + c\underline{D}_{EE}\{w(x,y)\}$$
$$= \underline{D}_{EE}\{g(x,y)\} + D_{EE}\{K(x,y)\}\underline{D}_{EE}\{w(x,y)\},$$

And

$$\sum_{h=1}^{l} a_h \overline{D}_{EE}\left\{\frac{\partial^h w(x,y)}{\partial x^h}\right\} + \sum_{j=1}^{m} b_j \overline{D}_{EE}\left\{\frac{\partial^j w(x,y)}{\partial y^j}\right\} + c\overline{D}_{EE}\{w(x,y)\}$$
$$= \overline{D}_{EE}\{g(x,y)\} + D_{EE}\{K(x,y)\}\overline{D}_{EE}\{w(x,y)\}.$$

Then from (6) and (7) we have

$$\left(\sum_{h=1}^{l} a_h (iv^n)^h + \sum_{j=1}^{m} b_j (iv^n)^j + c - D_{EE}\{K(x,y)\}\right) \underline{D}_{EE}\{w(x,y,\alpha)\}$$
$$= \underline{D}_{EE}\{g(x,y)\} + \sum_{h=1}^{l} \sum_{r=1}^{h} a_h (iv^n)^{h-r} \underline{E}^c\left\{\frac{\partial^{h-r} w}{\partial x^{h-r}}(0,y,\alpha)\right\}$$
$$+ \sum_{j=1}^{m} \sum_{r=1}^{j} b_j (iv^n)^{j-r} \underline{E}^c\left\{\frac{\partial^{j-r} w}{\partial y^{j-r}}(x,0,\alpha)\right\},$$

And

$$\left(\sum_{h=1}^{l} a_h (iv^n)^h + \sum_{j=1}^{m} b_j (iv^n)^j + c - D_{EE}\{K(x,y)\}\right) \overline{D}_{EE}\{w(x,y,\alpha)\}$$
$$= \underline{D}_{EE}\{g(x,y)\} + \sum_{h=1}^{l} \sum_{r=1}^{h} a_h (iv^n)^{h-r} \overline{E^c}\left\{\frac{\partial^{h-r} w}{\partial x^{h-r}}(0,y,\alpha)\right\}$$
$$+ \sum_{j=1}^{m} \sum_{r=1}^{j} b_j (iv^n)^{j-r} \overline{E^c}\left\{\frac{\partial^{j-r} w}{\partial y^{j-r}}(x,0,\alpha)\right\}.$$

Using the initial conditions we get

$$\left(\sum_{h=1}^{l} a_h (iv^n)^h + \sum_{j=1}^{m} b_j (iv^n)^j + c - D_{EE}\{K(x,y)\}\right) \underline{D}_{EE}\{w(x,y,\alpha)\}$$
$$= \underline{D}_{EE}\{g(x,y)\} + \sum_{h=1}^{l} \sum_{r=1}^{h} a_h (iv^n)^{h-r} \underline{E}^c\{\Phi_{h-r}(0,y,\alpha)\}$$
$$+ \sum_{j=1}^{m} \sum_{r=1}^{j} b_j (iv^n)^{j-r} \underline{E}^c\{\Psi_{j-r}(x,0,\alpha)\},$$

And

$$\left(\sum_{h=1}^{l} a_h (iv^n)^h + \sum_{j=1}^{m} b_j (iv^n)^j + c - D_{EE}\{K(x,y)\}\right)\overline{D}_{EE}\{w(x,y,\alpha)\}$$

$$= \underline{D}_{EE}\{g(x,y)\} + \sum_{h=1}^{l}\sum_{r=1}^{h} a_h (iv^n)^{h-r}\overline{E^c}\{\Phi_{h-r}(0,y,\alpha)\}$$

$$+ \sum_{j=1}^{m}\sum_{r=1}^{j} b_j (iv^n)^{j-r}\overline{E^c}\{\Psi_{j-r}(x,0,\alpha)\}.$$

Then

$$\underline{D}_{EE}\{w(x,y,\alpha)\}$$
$$= \frac{\underline{D}_{EE}\{g(x,y)\} + \sum_{h=1}^{l}\sum_{r=1}^{h} a_h(iu^n)^{h-r}\underline{E^c}\{\Phi_{h-r}(0,y,\alpha)\} + \sum_{j=1}^{m}\sum_{r=1}^{j} b_j(iv^n)^{j-r}\underline{E^c}\{\Psi_{j-r}(x,0,\alpha)\}}{\sum_{h=1}^{l} a_h (iu^n)^h + \sum_{j=1}^{m} b_j(iv^n)^j + c - D_{EE}\{K(x,y)\}}, \quad (9)$$

And

$$\overline{D}_{EE}\{w(x,y,\alpha)\}$$
$$= \frac{\overline{D}_{EE}\{g(x,y)\} + \sum_{h=1}^{l}\sum_{r=1}^{h} a_h(iu^n)^{h-r}\overline{E^c}\{\Phi_{h-r}(0,y,\alpha)\} + \sum_{j=1}^{m}\sum_{r=1}^{j} b_j(iv^n)^{j-r}\overline{E^c}\{\Psi_{j-r}(x,0,\alpha)\}}{\sum_{h=1}^{l} a_h (iu^n)^h + \sum_{j=1}^{m} b_j(iv^n)^j + c - D_{EE}\{K(x,y)\}}. \quad (10)$$

2. If $\widetilde{w}(x,y)$ be (ii)-differentiable, then

$$\underline{D}_{EE}\{w(x,y,\alpha)\}$$
$$= \frac{\underline{D}_{EE}\{g(x,y)\} + \sum_{h=1}^{l}\sum_{r=1}^{h} a_h(iu^n)^{h-r}\overline{E^c}\{\Phi_{h-r}(0,y,\alpha)\} + \sum_{j=1}^{m}\sum_{r=1}^{j} b_j(iv^n)^{j-r}\overline{E^c}\{\Psi_{j-r}(x,0,\alpha)\}}{\sum_{h=1}^{l} a_h (iu^n)^h + \sum_{j=1}^{m} b_j(iv^n)^j + c - D_{EE}\{K(x,y)\}},$$

And

$$\overline{D}_{EE}\{w(x,y,\alpha)\}$$
$$= \frac{\overline{D}_{EE}\{g(x,y)\} + \sum_{h=1}^{l}\sum_{r=1}^{h} a_h(iu^n)^{h-r}\underline{E^c}\{\Phi_{h-r}(0,y,\alpha)\} + \sum_{j=1}^{m}\sum_{r=1}^{j} b_j(iv^n)^{j-r}\underline{E^c}\{\Psi_{j-r}(x,0,\alpha)\}}{\sum_{h=1}^{l} a_h (iu^n)^h + \sum_{j=1}^{m} b_j(iv^n)^j + c - D_{EE}\{K(x,y)\}}.$$

Finally, we get $\widetilde{w}(x,y)$ by using the inverse of FDEET such that $\widetilde{w}(x,y) = \left(\underline{w}(x,y,\alpha), \overline{w}(x,y,\alpha)\right)$.

***Examples:***

In this section, we find the solution of partial convolution Volterra fuzzy integro-differential equation by using the previous method.

***Example1:***

Consider the following partial Volterra fuzzy integro-differential equation

$$\frac{\partial^2 \widetilde{w}(x,y)}{\partial x^2} \oplus \frac{\partial^2 \widetilde{w}(x,y)}{\partial y^2} \oplus \widetilde{w}(x,y) = \widetilde{g}(x,y) \oplus (F\mathcal{R})\int_0^y (F\mathcal{R})\int_0^x e^{x-\mu+y-\tau} \odot \widetilde{w}(\mu,\tau)\,d\mu\,d\tau, \quad (11)$$

$(x,y) \in [0,1] \times [0,1], \sigma \in [0,1]$.

With initial condition

$$\widetilde{w}(x,0) = \big(e^x(2+\sigma), e^x(4-r)\big), \quad \frac{\partial \widetilde{w}(x,0)}{\partial y} = \big(e^x(2+\sigma), e^x(4-r)\big),$$

$$\widetilde{w}(0,y) = \big(e^y(2+\sigma), e^y(4-r)\big), \quad \frac{\partial \widetilde{w}(0,y)}{\partial x} = \big(e^y(2+\sigma), e^y(4-r)\big),$$

And

$$\widetilde{g}(x,y) = \big(e^{x+y}(3-xy)(2+\sigma), e^{x+y}(3-xy)(4-r)\big).$$

In this case $l = m = 2, a_1 = b_1 = 0, a_2 = b_2 = c = 1,$

$K(x-\mu, y-\tau) = e^{x-\mu+y-\tau} > 0, \forall 0 \le \mu \le x \le 1 \ and \ 0 \le \tau \le y \le 1.$

Hence,

$$\Psi_0(x,0) = \big(e^x(2+\sigma), e^x(4-r)\big), \quad \Psi_1(x,0) = \big(e^x(2+\sigma), e^x(4-r)\big),$$

$$\Phi_0(0,y,\alpha) = \big(e^y(2+\sigma), e^y(4-r)\big), \quad \Phi_1(0,y) = \big(e^y(2+\sigma), e^y(4-r)\big). \ [27]$$

Apply double complex EE transform for both sides of eq. (9) and find the following

$$D_{EE}\{k(x,y)\} = D_{EE}\{e^{x+y}\} = \frac{(1+iu^n)(1+iv^n)}{(1+u^{2n})(1+v^{2n})},$$

$$\widetilde{D}_{EE}\{\Psi_0(x,0)\} = \widetilde{D}_{EE}\{\Psi_1(x,0)\} = \left(\frac{-(1+iu^n)}{(1+u^{2n})}(2+\sigma), \frac{-(1+iu^n)}{(1+u^{2n})}(4-r)\right),$$

$$\widetilde{D}_{EE}\{\Phi_0(0,y)\} = \widetilde{D}_{EE}\{\Phi_1(0,y)\} = \left(\frac{-(1+iv^n)}{(1+v^{2n})}(2+\sigma), \frac{-(1+iv^n)}{(1+v^{2n})}(4-r)\right).$$

From the theorem (5) we get:

$$\widetilde{D}_{EE}\{g(x,y)\} = \left(\frac{3}{(1-iu^n)(1-iv^n)} - \frac{1}{(1-iu^n)^2(1-iv^n)^2}(2+\sigma), \frac{3}{(1-iu^n)(1-iv^n)} - \frac{1}{(1-iu^n)^2(1-iv^n)^2}(4-r)\right).$$

Then by (9) and (10) for the solution of equation we have:

$$\underline{D}_{EE}\{w(x,y,\sigma)\} = \frac{(1+iu^n)(1+iv^n)}{(1+u^{2n})(1+v^{2n})}(2+\sigma), \overline{D}_{EE}\{w(x,y,\sigma)\} = \frac{(1+iu^n)(1+iv^n)}{(1+u^{2n})(1+v^{2n})}(4-r).$$

By using the inverse of double complex EE transform the solution of the equation will be

$$\widetilde{w}(x,y) = \big(e^{x+y}(2+\sigma), e^{x+y}(4-r)\big).$$

## *Conclusion*

In this paper, the concept of partial convolution Volterra fuzzy integro-differential equations have been introduced and solved by applying the double fuzzy complex EE transform. New results on DFEET for fuzzy partial H-derivative of the n-th order have been introduced.

By using parametric form of fuzzy functions we convert the investigated equation to a nonlinear partial Volterra integro-differential equation in a crisp case and applying DFEET for these equations to get algebraic equations. Hence we find lower and upper functions of the solution. Finally, the giving example shows that the investigation method is effective in solving the equations of considered kind.

*References:*

1. Chang, S.; Zadeh, L. On fuzzy mapping and control. IEEE Trans. Syst. Man Cybern. 1972, 1, 30–34.
2. Aslam, M.; Khan, N.; Ali Hussein, A.-M. Design of Variable Sampling Plan for Pareto Distribution Using Neutrosophic Statistical Interval Method. Symmetry 2019, 11, 80.
3. Aslam, M. A New Sampling Plan Using Neutrosophic Process Loss Consideration. Symmetry 2018, 10, 132.
4. Aslam, M.; Arif, O.H.; Khan Sherwan, A.R. New Diagnosis Test under the Neutrosophic Statistics :An Application to Diabetic Patients. BioMed Res. Int. 2020, 2020, 2086185.
5. Aslam, M. Attribute Control Chart Using the Repetitive Sampling Under Neutrosophic System. IEEE Access, 2019، 7، 15374–15367.
6. Allahviranloo, T.; Kermani, M. Numerical methods for fuzzy linear partial differential equations under new definition for derivative. Iran. J. Fuzzy Syst. 2010, 7, 33–50.
7. Chalco-Cano, Y.; Roman-Flores, H. On new solutions of fuzzy differential quations. Chaos Solitons Fractals, 2008، 38، 112–119.
8. Alidema, A.; Georgieva, A. Adomian decomposition method for solving two-dimensional nonlinear Volterra fuzzy integral equations. AIP Conf. Proc. 2018, 2048, 050009.
9. Georgieva, A. Solving two-dimensional nonlinear Volterra-Fredholm fuzzy integral equations by using Adomian decomposition method. Dyn. Syst. Appl. 2018, 27, 819–837.
10. Abbasbandy, S.; Hashemi, M. A series solution of fuzzy integro-differential equations. J. Fuzzy Set Valued Anal. 2012, 2012, jfsva-00066.
11. Hooshangian, L. Nonlinear Fuzzy Volterra Integro-differential Equation of N-th Order: Analytic Solution and Existence and Uniqueness of Solution. Int. J. Ind. Math. 2019, 11, 12.
12. Ahmad, J.; Nosher, H. Solution of Different Types of Fuzzy Integro-Differential Equations Via Laplace Homotopy Perturbation Method. J. Sci. Arts 2017, 17, 5.
13. Mikaeilv, N.; Khakrangin, S.; Allahviranloo, T. Solving fuzzy Volterra integro-differential equation by fuzzy differential transform method. In Proceedings of the 7th Conference of the European Society for Fuzzy Logic and Technology, Aix-Les-Bains, France, 18–22 July 2011; pp. 18–22.
14. Arnoldus, H.F. Application of the magnetic field integral equation to diffraction and reection by a conducting sheet. Int. J. Theor. Phys. Group Theory Nonlinear Opt. 2011, 14, 1–12.
15. Ma, S.Q.; Chem, F.C.; Zhao, Z.Q. Choquet type fuzzy complex-values integral and its application in classification. Fuzzy Eng. Oper. Res. 2012, 147, 229-–237.
16. Debnath L., and Bhatta D. , "Integral Transforms and Their Applications", 2nd ed. Chapman & Hall/CRC, 2007.
17. F. Bin Muhammed Belgacem, "Sumudu transform fundamental properties investigations and applications," Journal of Applied Mathematics and Stochastic Analysis, vol. 2006, pp. 1-23, 2006, 10.1155/JAMSA/2006/91083.

18 Ali Hassan Mohammed., Athera Nemakathem, Solving Euler's Equation by using New Transformation, karbala University magazine for Completely sciences, Vol.4, N0.2, 2008, pp. 103-109.
19 Mohammed A. H., Makttoof S. F., A Complex Al-Tememe Transform, International J. of Pure & Engg. Mathematics, Vol. 5, No. II , 2017, pp. 17-30.
20 Mohand M., Abdelrahim Mahgoub, The New Integral Transform "Mahgoub Transform", Advances in Theoretical and Applied Mathematics, Vol.11, No.4, 2016, pp. 391- 398.
21 Mohand M. Abdelrahim Mahgoub, The New Integral Transform "Mohand Transform", Advances in Theoretical and Applied Mathematics, Vol. 12, No.2, 2017, pp.113- 120.
22 Eman A. Mansour, Sadiq Mehdi, Emad A. Kuffi, The new integral transform and Its Applications, Int. J. Nonlinear Anal. Appl. Vol.12, No.2, 2021, pp. 849-856.
23 Eman A. Mansour, Emad A. Kuffi, Sadiq A. Mehdi, The new integral transform "SEE transform " and its Applications, Periodicals of Engineering and Natural Sciences, Vol.9, N0.2, 2021, pp.1016-1029.
24 E. A. Kuffi, Abbas E.s., A complex integral transform "Complex EE transform" and its application on Differential Equations, Mathematical Statistician and Engineering Applications, vol. 71, No.2, pp. 263-266, 2022.
25 Alidema, A.; Georgieva, A. Applications of the double fuzzy Sumudu transform for solving Volterra fuzzy integral equations. AIP Conf. Proc. 2019, 2172, 060006 .
26 Mohand, M.; Mahgob A. Solution of Partial Integro-Differential Equations by Double Elzaki Transform Method. Math. Theory Model. 2015, 5, 61–65.
27 Atanaska Georgieva , Double Fuzzy Sumudu Transform to Solve Partial Volterra Fuzzy Integro-Differential Equations,Mathematics 2020, 8, 692, doi:10.3390/math8050692.
28 Ahmad Issa, E. A. Kuffi, Abbas, On The Double Integral Transform (Complex EE Transform) and Their Properties and Applications, Ibn Al-Haitham Journal for Pure and Applied Sciences, 2023.
29 Abaas S. T., Alkiffai A. N., Kadhim A. N., Efficient integral fuzzy transformation for practical problems, A thesis of Ph.D. submitted to the council of University of Kufa, Faculty of Education for Girls, 2022.